\documentclass [twoside,12pt,draft] {article}

\usepackage{amsfonts}
\usepackage{amssymb}
\usepackage[arrow, matrix, curve]{xy}

\newtheorem{theorem}{\\ $\ \ \ \ $ Theorem}[section]
\newtheorem{lemma}[theorem]{\\ $\ \ \ \ $ Lemma}
\newtheorem{definition}[theorem]{\\ $\ \ \ \ $ Definition}
\newtheorem{remark}[theorem]{\\ $\ \ \ \ $ Remark}
\newtheorem{example}[theorem]{\\ $\ \ \ \ $ Example}

\newcommand{\adj}{\mathrm{adj}}
\newcommand{\C}{\mathbb{C}}
\newcommand{\cl}{\mathrm{cl}}
\newcommand{\iD}{\mathit{\Delta}}

\newcommand{\Id}{\mathit{Id}}
\newcommand{\iL}{\mathit{\Lambda}}

\newcommand{\iO}{\mathit{\Omega}}
\newcommand{\iP}{\mathit{\Pi}}
\newcommand{\PsDO}{\mathit{\Psi}}
\newcommand{\R}{\mathbb{R}}

\begin{document}

\title{Elliptic Quasicomplexes 
\\
       on Compact Closed Manifolds}

\date{\today}

\author{D. Wallenta}

\maketitle

\begin{abstract}
We consider quasicomplexes of pseudodifferential operators on a smooth compact manifold without boundary. 
To each quasicomplex we associate a complex of symbols. 
The quasicomplex is elliptic if this symbol complex is exact away from the zero section. 
We prove that elliptic quasicomplexes are Fredholm. 
Moreover, we introduce the Euler characteristic for elliptic quasicomplexes and prove a generalization of the Atiyah-Singer index theorem.
\end{abstract}



\section{Introduction}
\label{s.i}
\setcounter{equation}{0}

In this article the theory of elliptic quasicomplexes on compact closed manifolds is developed. 
This is a generalization of the theory of elliptic complexes, which are studied since the 1950s by A. Grothendieck,
                   D.~C. Spencer,
                   M.~F. Atiyah, 
                   R. Bott, 
and others. 
The standard example of such a complex is the de Rham complex.

Complexes, also called cochain complexes, are well-known mathematical objects. In the first two sections we will recall basic facts about general Fredholm complexes in Banach and Hilbert spaces, respectively. 

One may ask what happens with complexes under ``small'' perturbations of their differentials.  
This leads us to the theory of quasicomplexes introduced by Putinar in  
   \cite{Puti82} 
and developed by Tarkhanov at the end of the 1990s. 
In Section \ref{s.Quasicomplexe} we present a summary of the modern theory of Fredholm quasicomplexes in Hilbert spaces. 
An operator is said to be ``small'' in this context, if it is compact. 
In this section the main results of	\cite{Puti82}, 
                                    \cite{Tark07} and 
                                    \cite{KrupTarkTuom07} 
are discussed. 
Note that it depends on the structure of the underlying spaces whether or not an operator is ``small.''  

The concept of classical pseudodifferential operators on a compact closed manifold is sketched in Section \ref{s.Pseudo}. 
Such operators are very simple to manage, see \cite{Shub87}. 
For example, they form an algebra with unit and have principal symbols living on the cotangent bundle of the manifold, like linear partial differential operators. 
The ellipticity of a pseudodifferential operator is defined via its symbol. Moreover, each of these operators extends by continuity to a linear operator in appropriate Sobolev spaces. 
   
Complexes of pseudodifferential operators operating on spaces of smooth sections of vector bundles are introduced in Section \ref{s.Index}. 
Note that these spaces are (Fr\'{e}chet-Schwartz spaces and so) no longer Banach, hence compact operators fail to be ``small'' in this context.
Following \cite{KrupTarkTuom07} we call a pseudodifferential operator ``small'' if it is smoothing. 
According to this, we elaborate the theory of elliptic quasicomplexes in Section \ref{s.Elliptisch}. 
We define quasicomplexes and their parametrizes.
We prove that the concept of ellipticity can be extended naturally to quasicomplexes of pseudodifferential operators. 
Moreover, we introduce the Euler characteristic of an elliptic quasicomplex and show a generalization of the Atiyah-Singer index formula.

Note that each elliptic operator gives rise to an elliptic (quasi-) complex of length $1$. 
Since many important operators of theoretical physics, for example the Laplacian, are elliptic, one can find practical applications in this field,
   see for instance \cite{Schw83}.

\section{Complexes}
\label{s.Complexes}
\setcounter{equation}{0}

Complexes of operators are generalizations of single operators. 
If $A \! : V \to W$ is a linear map between vector spaces, then $A$ defines the so-called short complex
$$  
   0  
 \rightarrow 
   V 
 \stackrel{A}{\rightarrow}  
   W 
 \rightarrow 
   0.
$$

By a (cochain) complex $V^{\cdot}$ is meant a sequence of linear maps between vector spaces  
$$ 
   V^{\cdot}\!:\,\, 
   0 
 \rightarrow  
   V^{0} 
 \stackrel{A^{0}}{\rightarrow}  
   V^{1} 
 \stackrel{A^{1}}{\rightarrow}  
   \ldots
 \stackrel{A^{N-1}}{\rightarrow}  
   V^{N} 
 \rightarrow
   0
$$
with $A^{i+1} A^{i}=0$ for all $i = 0, 1, \ldots, N-1$.
For such a complex $V^{\cdot}$, we set 
   $V^i = 0$ for $i \in \mathbb{Z} \setminus \{0,...,N\}$ and 
   $A^i = 0$ for $i \in \mathbb{Z} \setminus \{0,...,N-1\}$. 
To each complex the differential $A$ is associated by 
   $A v = A^{i} v$ 
for $v \in V^{i}$. 
Since $A^{2} = 0$ the differential is nilpotent.   
We will write $(V^{\cdot},A)$ instead of $V^\cdot$, if we want to emphasize which differential is used. 

For any $v \in V^{i-1}$, 
$$
   A^{i} (A^{i-1} v) 
 = (A^{i} A^{i-1}) v
 = 0
$$  
whence
$
   \mathrm{im}\, A^{i-1} \subset \mathrm{ker}\, A^{i}, 
$
i.e., the image of $A^{i-1}$ is a subspace of the kernel of $A^{i}$.  
The quotient space
$$
   H^{i} (V^{\cdot})
 :=
   \mathrm{ker}\, A^{i} / \mathrm{im}\, A^{i-1}
$$
is called the cohomology of the complex at step $i$. 
A complex is said to be exact at step $i$ if
   $H^{i} (V^{\cdot}) = 0$.

\begin{remark}
\label{B.kurzerKomplex}
{\em 
For a short complex we have
$$
\begin{array}{ccccccc}
   H^{0}(V^{\cdot})
 & = 
 & \mathrm{ker}\, A^{0} / \mathrm{im}\, A^{-1}
 & =
 & \mathrm{ker}\, A / \{0\} 
 & \cong 
 & \mathrm{ker}\, A,
\\
   H^{1} (V^{\cdot})
 & =
 & \mathrm{ker}\, A^{1} / \mathrm{im}\, A^{0}
 & =
 & W / \mathrm{im}\, A
 & =:
 & \mathrm{coker}\, A.
\end{array}
$$
}
\end{remark}

\begin{definition}
Let $V^{\cdot}$ be a complex with finite dimensional cohomology. 
We define the Euler characteristic of the complex by 
$$
   \chi (V^{\cdot})
 :=
   \sum_{i} (-1)^{i}\, \mathrm{dim}\, H^{i} (V^{\cdot}).
$$
\end{definition}
   
\begin{example}
{\em 
If $V^{\cdot}$ is a short complex, then from Remark \ref{B.kurzerKomplex} it follows that 
$$  
   \chi (V^{\cdot})
 = \mathrm{dim}\, H^{0} - \mathrm{dim}\, H^{1} 
 = \mathrm{dim}\, \ker A - \mathrm{dim}\, \mathrm{coker}\, A,
$$
which is the index of the operator $A$.
}
\end{example}

Let $(V^{\cdot},A)$ and
    $(W^{\cdot},B)$ 
be two complexes. 
Without restriction of generality we can take complexes of the same length $N$.    
A homomorphism of the complexes $V^{\cdot}$ and 
                                $W^{\cdot}$ 
is a sequence of linear maps 
    $L^{i}\!: V^{i} \to W^{i}$
which makes the diagram  
$$
\begin{array}{ccccccccccc}
   {\displaystyle 0}
 & \rightarrow
 & {\displaystyle V^0}
 & \stackrel{A^{0}}\rightarrow
 & {\displaystyle V^1}
 & \stackrel{A^{1}}\rightarrow
 & {\displaystyle \ldots}
 & \stackrel{A^{N-1}}\rightarrow
 & {\displaystyle V^N}
 & \rightarrow
 & {\displaystyle 0}
\\
 & 
 & \mbox{ }\downarrow L^0 
 & 
 & \mbox{ }\downarrow L^1
 & 
 & 
 & 
 & \mbox{ }\downarrow L^N
 & 
 & 
\\
   {\displaystyle 0}
 & \rightarrow
 & {\displaystyle W^0}
 & \stackrel{B^{0}}\rightarrow
 & {\displaystyle W^1}
 & \stackrel{B^{1}}\rightarrow
 & {\displaystyle \ldots}
 & \stackrel{B^{N-1}}\rightarrow
 & {\displaystyle W^N}
 & \rightarrow
 & {\displaystyle 0}
\end{array}
$$
commutative, i.e.
$
   L^{i+1} A^{i} = B^{i} L^{i}
$ 
for all $i = 0, 1, \ldots, N-1$. 
Each homomorphism $\{ L^{i} \}$ induces a sequence of homomorphisms 
$
   HL^{i}\!: H^{i} (V^{\cdot}) \to H^{i} (W^{\cdot}) 
$ 
of the cohomology by 
$
   HL^{i} [v]: = [L^{i} v]
$
for $[v] \in H^{i} (V^{\cdot})$.
It is easy to see that this is well defined.   
The homomorphisms of $V^{\cdot}$ to 
                     $V^{\cdot}$ 
itself are called endomorphisms of this complex.

Suppose $V^{\cdot}$ is a complex with finite dimensional cohomology and  $\{E^{i}\}$ an endomorphism of the complex. 
Then $HE^{i}$ is an endomorphism of the finite dimensional space 
$H^{i} (V^{\cdot})$, and so the trace $\mathrm{tr}\, HE^{i}$ is well defined.  
The alternating sum 
$$ 
   L (E) := \sum_{i} (-1)^{i}\, \mathrm{tr}\, HE^{i}
$$ 
is called the Lefschetz number of the endomorphism.
If $E^i = \Id_{V^i}$ are the identity maps, then the trace  
   $\mathrm{tr}\, HE^{i}$
just amounts to the dimension of $H^{i} (V^{\cdot})$ whence 
$ 
   L (\Id_{V^\cdot}) = \chi (V^{\cdot}).
$
 
We have established complexes as sequences of linear maps between vector spaces which is adequate for algebraic analysis. 
If we want to use methods of calculus, we have to include continuous linear maps between topological vector spaces. 
For the rest of the paper we will understand complexes in this topological sense. 
If we use homomorphisms  $\{ L^i \}$, we will suppose the maps $L^i$ to be continuous, too.      
 
\begin{example}
{\em 
Let $X$ be a smooth (i.e. $C^\infty$) manifold of dimension $n$. 
Denote by $\iO^{q} (X)$ the space of differential forms of degree $q$ with
smooth coefficients on $X$ and 
   $d \! : \iO^{q} (X) \to \iO^{q+1}(X)$ 
the exterior derivative.  
Locally any $\omega \in \iO^{q} (X)$ looks like 
$$
   \omega (x)
 = \sum_{J = (j_1, \ldots, j_q) \atop 
         1 \leq j_1 < \ldots < j_q \leq n}  
   \omega_{J} (x)\, dx^J 
$$
for $x = (x^1, \ldots, x^n)$ in a coordinate patch $U$ of $X$, where 
   $dx^J = dx^{j_1} \wedge \ldots \wedge dx^{j_q}$ and  
   $\omega_{J} \in C^\infty (U,\mathbb{R})$.
The derivative is given by 
$$
   d \omega (x)
 = \sum_{J = (j_1, \ldots, j_q) \atop 
         1 \leq j_1 < \ldots < j_q \leq n}  
   d \omega_{J} (x) \wedge dx^J 
$$ 
for $x \in U$.
It is linear and satisfies $d^{2}=0$. 
Hence 
$$
   \iO^{\cdot} (X)\!:\,\, 
   0 
 \rightarrow 
   \iO^{0}(X) 
 \stackrel{d}{\rightarrow}   
   \iO^{1}(X)
 \stackrel{d}{\rightarrow}  
   \ldots
 \stackrel{d}{\rightarrow}    
   \iO^{n} (X) 
 \rightarrow 
   0
$$
is a complex.
This complex is referred to as the de Rham complex of $X$ and its cohomology 
$
   H^{i}_{\mathit{dR}} (X)  
 :=  
   H^{i} (\iO^{\cdot} (X))
$
are called the de Rham cohomology of $X$.
}
\end{example}

The de Rham complex is a classical example of complexes.
The numbers 
   $\mathrm{dim}\, H^{i}_{\mathit{dR}} (X)$
are called Betti numbers of the underlying manifold $X$.
They depend on certain topological properties of $X$. 
We will come back to this example later.

\section{Fredholm Property}
\label{Fred}
\setcounter{equation}{0}

In this section we mainly use functional analytic methods. 
First we discuss Fredholm operators. 

Let $V$ and $W$ be Banach spaces. 
An operator $A \in \mathcal{L} (V,W)$ is called Fredholm if both
   $\mathrm{ker}\, A$ and  
   $\mathrm{coker}\, A$ 
are of finite dimension. 
The index of a Fredholm operator $A$ is the number
$ 
\mathrm{ind}\, A := \mathrm{dim}\, \ker A -
                    \mathrm{dim}\, \mathrm{coker}\, A.
$
As usual, we write $\mathcal{K} (V,W)$ for the set of all compact operators acting from $V$ to $W$.
Note that this is a closed subspace of the Banach space $\mathcal{L} (V,W)$ and hence a Banach space itself. 
The composition of a compact operator and a bounded linear operator is always compact. 
In particular $\mathcal{K} (V):=\mathcal{K} (V,V)$ is an ideal in 
              $\mathcal{L} (V):=\mathcal{L} (V,V)$.

\begin{theorem}
Let $A \in \mathcal{L} (V,W)$ be Fredholm and 
    $K \in \mathcal{K} (V,W)$. 
Then $A + K$ is Fredholm and
$
   \mathrm{ind}\, A = \mathrm{ind}\, (A+K).
$ 
\end{theorem}

An operator $P \in \mathcal{L} (W,V)$ is called a parametrix of 
            $A \in \mathcal{L} (V,W)$, 
if 
$$
\begin{array}{rcl}
   \Id_V - P A
 & \in 
 & \mathcal{K} (V),
\\
   \Id_W - A P
 & \in 
 & \mathcal{K} (W)
\end{array}
$$ 
holds. 
In other words, by a parametrix of $A$ is meant an inverse of $A$ modulo compact operators. 
This property can be described by using a familiar construction with quotient spaces which goes back as far as \cite{Calk41}. 
Given a Banach space  $\Sigma$, we set 
$$
\begin{array}{rcl}
   \phi_{\Sigma} (V)
 & := 
 & \mathcal{L} (\Sigma,V) / \mathcal{K} (\Sigma,V),
\\
   \phi_{\Sigma} (W)
 & :=
 & \mathcal{L} (\Sigma,W) / \mathcal{K} (\Sigma,W).
\end{array}
$$
Furthermore, for $A \in \mathcal{L} (V,W)$, we introduce a map 
$
   \phi_{\Sigma} (A)\!: \phi_{\Sigma} (V) \to \phi_{\Sigma} (W)
$    
by
$$
   \phi_{\Sigma} (A) [O] := [A \circ O]
$$
for $O \in \mathcal{L} (\Sigma,V)$. 
This defines a functor $\phi_{\Sigma}$ 
   from the category of Banach spaces 
   to the category of `Banach algebras', 
such that

\medskip

   $i)$  
$\phi_{\Sigma} (A) = 0$, if $A$ is compact;

   $ii)$ 
$\phi_{\Sigma} (BA) = \phi_{\Sigma} (B)\, \phi_{\Sigma}(A)$;

   $iii)$ 
$\phi_{\Sigma} (\Id_V) = \Id_{\phi_{\Sigma} (V)}$

\medskip
\noindent
for all
   $A \in \mathcal{L} (V,W)$ and
   $B \in \mathcal{L} (W,Z)$. 
If $\Sigma = V$ then 
   $\phi_{\Sigma} (V) = \mathcal{L} (V) / \mathcal{K} (V)$ 
is a Banach algebra indeed.

Let $A \in \mathcal{L} (V,W)$.
The operator $\phi_{\Sigma} (A)$ proves to be invertible for each Banach space $\Sigma$ if and only if there is an operator 
   $P \in \mathcal{L} (W,V)$ 
with the property that
$$
\begin{array}{rcl}  
   \phi_{\Sigma} (P)\, \phi_{\Sigma} (A)
 & =
 & \Id_{\phi_{\Sigma} (V)},
\\
   \phi_{\Sigma} (A)\, \phi_{\Sigma} (P) 
 & =
 & \Id_{\phi_{\Sigma} (W)}
\\
\end{array}
$$
for all Banach spaces $\Sigma$. 

\begin{theorem}
\label{S.FredOp}
Let $V$ and $W$ be Banach spaces and $A \in \mathcal{L} (V,W)$. 
The following are equivalent:

\medskip 

   $i)$  
A is Fredholm.

   $ii)$ 
$A$ possesses a parametrix.
  
	 $iii)$ 
$\phi_{\Sigma} (A)$ is invertible for each Banach space $\Sigma$.
\end{theorem}

The equivalence of i) and ii) is known as theorem of Th. Atkinson. 
It is easy to see that an operator $A \in \mathcal{L} (V,W)$  is Fredholm if and only if the short complex    
$$
   V^{\cdot} :\,\, 
   \,0\, 
 \rightarrow 
   \,V\, 
 \stackrel{A}{\rightarrow}  
   \,W\,  
 \rightarrow 
   \,0
$$
possesses finite dimensional cohomology. 
Now, we are in a position to extend the concept of Fredholm operators to complexes of Banach spaces in a natural way.

\begin{definition}
A complex $V^{\cdot}$ is said to be Fredholm if its cohomology
   $H^{i} (V^\cdot)$ 
is finite dimensional at each step $i$.
\end{definition}
 
Let $(V^{\cdot}, A)$ be a complex of Banach spaces and continuous linear maps. 
By a parametrix of this complex is meant a sequence of operators 
   $P^{i} \in \mathcal{L} ( V^{i}, V^{i-1})$ 
satisfying 
$$ 
   P^{i+1} A^{i} + A^{i-1} P^{i} = \Id_{V^{i}} - K^{i}
$$ 
with 
   $K^{i} \in \mathcal{K} (V^{i})$ 
for all $i = 0, 1, \ldots, N$.

Obviously, each Fredholm complex $V^{\cdot}$ possesses an Euler characteristic. 
Moreover, the Lefschetz number $L (E)$ is well defined for any endomorphism $E$ of $V^{\cdot}$. 
Theorem \ref{S.FredOp} says that a short complex is Fredholm if and only if this complex possesses a parametrix. 
This result extends naturally to arbitrary complexes of Banach spaces.
The proof is especially elegant for complexes of Hilbert spaces $V$, where one uses the so-called adjoint complex.
By this is meant the complex  
$$
   V^{\cdot}{}^\ast :\,\,
   \,0\,
 \leftarrow
   \,V^0\,
 \stackrel{A^0{}^\ast}{\leftarrow}
   \,V^1\,
 \stackrel{A^1{}^\ast}{\leftarrow}
   \,\ldots\,
 \stackrel{A^{N-1}{}^\ast}{\leftarrow}
   \,V^N\,
 \leftarrow
   \,0,
$$ 
where $A^i{}^\ast \in \mathcal{L} (V^{i+1},V^i)$ stands for the adjoint of $A^i$ in the sense of Hilbert spaces.
The equality $(A^\ast)^2 = 0$ is clear from $A^2 = 0$.
After W. V. D. Hodge, the operators  
$
   \iD^i = A^{i-1} A^{i-1}{}^\ast + A^i{}^\ast A^i
$ 
are called the Laplacians of the complex. 
The null-space of $\iD^i$ consists of all $h \in V^i$ satisfying 
   $A^i h = 0$ and
   $A^{i-1}{}^\ast h = 0$,
as is easy to see.

\begin{theorem}
\label{l.Laplacian}
Let $V^{\cdot}$ be a complex of continuous linear operators between Hilbert spaces. 
The following are equivalent: 

   $i)$  
$V^\cdot$ is Fredholm.

   $ii)$ 
All Laplacians $\iD^i$ of $V^\cdot$ are Fredholm.
 
	$iii)$ 
$V^\cdot$ possesses a parametrix.
\end{theorem}

The natural map $\ker \iD^i \mapsto H^i (V^\cdot)$ is obviously injective.
If $V^\cdot$ is Fredholm then this map is actually surjective, for the strong orthogonal decomposition 
$
   V^i = \ker \iD^i \oplus \mathrm{im}\, \iD^i
$
holds.
Indeed, denote by $H^i$ the orthogonal projection of $V^i$ onto the null-space of $\iD^i$.
On introducing the so-called Green operator $G^i$ in $V^i$ by
 $$
   G^i 
 := 
   (\iD^i \restriction_{(\ker \iD^i)^\perp})^{-1}\,  (\Id_{V^i} - H^i)
$$
we get
\begin{equation}
\label{eq.Hodge}
   \Id_{V^i}
 = H^i 
 + A^{i-1}\, (A^{i-1}{}^\ast G^i) 
 + (A^i{}^\ast G^{i+1})\, A^i
\end{equation}
for all $i = 0, 1, \ldots, N$,
the summands on the right-hand side being orthogonal.
This fact is known as abstract Hodge theory.

\section{Quasicomplexes}
\label{s.Quasicomplexe}
\setcounter{equation}{0}

Let $(V^\cdot,A)$ be a complex. 
We want to perturb the operators by ``small'' operators.
In the general setting of Banach spaces compact operators are ``small". 
We may perturbate the differential $A$ by a compact operator $K$. 
Formally we find  
$$
   (A + K)(A + K)
 = A^2 + AK + KA + K^2
 = AK + KA + K^2.
$$
Hence, the product is a compact operator, too. 

\begin{definition}
A sequence $(V^{\cdot},A)$ of operators
   $A^i \in \mathcal{L} (V^i,V^{i+1})$
in Banach spaces is called quasicomplex if its curvature is ``small'' at each step, i.e. 
$
   A^{i+1} A^i \in \mathcal{K} (V^i,V^{i+2})
$ 
for all $i = 0, 1, \ldots, N-2$.
\end{definition}

By the above, we obtain quasicomplexes by perturbing the differential of a complex by compact operators. 
In particular, each (finite) sequence of linear maps between finite-dimensional vector spaces is a quasicomplex.   

For quasicomplexes the cohomology is no longer defined, since the image of $A^{i-1}$ fails in general to lie in the null-space of $A^i$. 
However, in order to define Fredholm quasicomplexes we may use the construction with Calkin quotient spaces of Section \ref{Fred}. 
We choose an arbitrary Banach space $\Sigma$ and consider the sequence
$$
   \phi_{\Sigma} (V^{\cdot}) :\,\, 
   \,0\, 
 \rightarrow 
   \,\phi_{\Sigma} (V^{0})\, 
 \stackrel{\phi_{\Sigma} (A)}{\rightarrow} 
   \,\phi_{\Sigma} (V^{1})\, 
 \stackrel{\phi_{\Sigma}(A)}{\rightarrow}  
   \,\ldots\, 
 \stackrel{\phi_{\Sigma} (A)}{\rightarrow}  
   \,\phi_{\Sigma} (V^{N})\, 
 \rightarrow
   \,0
$$
which actually proves to be a complex.
Indeed, since $A^2$ is compact, we get 
$
   \phi_{\Sigma} (A) \phi_{\Sigma} (A)
 = \phi_{\Sigma} (A^2)
 = 0.
$  

A quasicomplex $V^{\cdot}$ is called Fredholm if 
   $\phi_{\Sigma} (V^{\cdot})$ 
is exact for each Banach space $\Sigma$.
The parametrix of a quasicomplex is defined in the same way as that of a complex.
A quasicomplex of Hilbert spaces is Fredholm if and only if it has a parametrix,
   see \cite{Tark07}.

\begin{theorem}
When perturbing the operators of a Fredholm complex of Hilbert spaces $(V^\cdot,A)$ by compact operators $C^i$, we obtain a Fredholm quasicomplex.
\end{theorem}

{\bf Proof.} 
Suppose $\{ P^i \}$ is a parametrix of $(V^\cdot,A)$, i.e.
   there are compact operators $K^i \in \mathcal{K} (V^i)$, such that
$$
   P^{i+1} A^i + A^{i-1} P^i = \Id_{V^i} - K^i
$$ 
for all $i = 0, 1, \ldots, N$.
We set $B^i := A^i + C^i$ and conclude that 
\begin{eqnarray*}
   P^{i+1} B^i + B^{i-1} P^i 
 & = &
   P^{i+1} A^i + A^{i-1} P^i + P^{i+1} C^i + C^{i-1} P^i
\\
 & = &
   \Id_{V^i}
 - (K^i - P^{i+1} C^i - C^{i-1} P^i),
\end{eqnarray*}
the operators in the parentheses on the right-hand side being in
   $\mathcal{K} (V^i)$.
Hence, $\{ P^i \}$ is a parametrix of the quasicomplex $(V^\cdot,B)$. 
This implies that $(V^\cdot,B)$ is Fredholm.
\hfill $\square$

The inverse theorem is also true but its proof is much more difficult.
The proof uses tricky the abstract Hodge decomposition (\ref{eq.Hodge}).
   
\begin{theorem}
\label{S.Reduktion}
Let $(V^{\cdot},A)$ be a Fredholm quasicomplex of Hilbert spaces.
Then, there exist operators  
   $D^{i} \in \mathcal{L} (V^{i},V^{i+1})$, 
such that $D^{i} - A^i \in \mathcal{K} (V^{i},V^{i+1})$ and 
          $D^{i+1} D^i = 0$.
\end{theorem}

This is the main result of \cite{Tark07}. 
The reduced complex may be used to define the Euler characteristic of a Fredholm quasicomplex. 
Let $(V^{\cdot},A)$ be a Fredholm quasicomplex of Hilbert spaces.
One introduces the Euler characteristic of $(V^{\cdot},A)$ by
\begin{equation}
\label{eq.Euler}
   \chi (V^{\cdot},A) := \chi (V^{\cdot},D),
\end{equation}
where $(V^{\cdot},D)$ is a reduced complex of $(V^{\cdot},A)$.
This definition is correct, i.e. the right-hand side of (\ref{eq.Euler}) doesn't depend on the particular choice of reduced complex.

The Hodge theory for Fredholm quasicomplexes was first elaborated in the paper \cite{KrupTarkTuom07}. 
Given a quasicomplex $(V^{\cdot},A)$ with Hilbert spaces $V^i$, the selfadjoint operators
$
   \iD^i = A^{i-1} A^{i-1}{}^\ast + A^i{}^\ast A^i  
$ 
are called the Laplacians of $(V^{\cdot},A)$.  

\begin{theorem}
\label{S.Fred Quasi}
Let $V^{\cdot}$ be a quasicomplex of Hilbert spaces. 
The following are equivalent: 

   $i)$  
$V^\cdot$ is Fredholm.

   $ii)$
The Laplacians $\iD^i$ of $V^\cdot$ are Fredholm operators. 
\end{theorem}

{\bf Proof.}
$i) \Rightarrow ii)$ 
Suppose $(V^{\cdot},A)$ is Fredholm.
By Theorem \ref{S.Reduktion} one can reduce this quasicomplex to a Fredholm complex $(V^{\cdot},D)$. 
The Laplacians of $(V^{\cdot},D)$ are Fredholm operators, which is due to
   Theorem \ref{l.Laplacian}.
Since the Laplacians of $(V^{\cdot},A)$ differ from those of
                        $(V^{\cdot},D)$
by compact operators, they are Fredholm, too. 

$ii) \Rightarrow i)$ 
Let the Laplacians $\iD^i$ be Fredholm. 
We choose a parametrix $G^i$ for $\iD^i$, i.e. a continuous linear operator in $V^i$ with
$$
\begin{array}{rcl}
   G^i \iD^i & = & \mathit{Id}_{V^i} - K^i_{l},
\\
   \iD^i G^i & = & \mathit{Id}_{V^i} - K^i_{r},
\end{array}
$$
where $K^i_{l}$ and $K^i_{r}$ are compact operators in $V^i$.
We know that $A^i \iD^i = \iD^{i+1} A^i$ holds modulo compact operators.
Multiplying this equality by $G^{i+1}$ from the left and
                          by $G^i$ from the right
we obtain 
$$
   G^{i+1} A^i = A^i G^i
$$
modulo compact operators from $V^i$ to $V^{i+1}$.
We set $P^i = G^{i-1} A^{i-1}{}^\ast$ and find
\begin{eqnarray*}
   P^{i+1} A^i + A^{i-1} P^i 
 & = &
   G^{i} A^{i}{}^\ast A^i + A^{i-1} G^{i-1} A^{i-1}{}^\ast
\\
 & = &
   G^{i} A^{i}{}^\ast A^i + G^{i} A^{i-1} A^{i-1}{}^\ast
\\
 & = &
	G^i \iD^i 
\\
 & = &
  \Id_{V^i}  
\end{eqnarray*}
modulo compact operators in $V^i$, i.e. the family
   $P^i : V^i \to V^{i-1}$ 
is a parametrix of $V^\cdot$.
Hence $V^\cdot$ is Fredholm. 
\hfill $\square$

\medskip

Similar to the Hodge theory for Fredholm complexes we 
   denote by $H^i$ the orthogonal projection of $V^i$ onto the null-space of  
   $\iD^i$ 
and 
   introduce the Green operator
$$
   G^i 
 := 
   (\iD^i \restriction_{(\ker \iD^i)^\perp})^{-1} (\Id_{V^i} - H^i).
$$
Then
$
   \Id_{V^i} = H^i + \iD^i G^i,
$
the summands on the right are orthogonal,
   cf. (\ref{eq.Hodge}).
It is easy to see that the operators $G^i$ satisfy 
   $A^i G^i - G^{i+1} A^i \in \mathcal{K} (V^i,V^{i+1})$. 
Hence it follows that the operators 
   $P^i := A^{i-1}{}^\ast G^i$ 
define a parametrix for $(V^{\cdot},A)$.

The theory of quasicomplexes would gain in interest if we introduce reasonably the Lefschetz number for endomorphism of quasicomplexes.
This question is at present far from being solved.
To highlight the problem, we consider a Fredholm quasicomplex $(V^\cdot,A)$ of Hilbert spaces and an endomorphism $\{ E^i \}$, i.e.
$
   E^{i+1} A^{i} = A^{i}E^{i}.
$
The idea is now to use any reduced complex $(V^\cdot,D)$ for $(V^\cdot,A)$ to define the Lefschetz number. 
Since $D^i = A^i + C^i$ implies 
$
   E^{i+1} D^{i} = D^{i} E^{i}
$
modulo $\mathcal{K} (V^i,V^{i+1})$, the sequence $\{ E^i \}$ fails to determine an endomorphism of $(V^\cdot,D)$. 
So nothing changes if we deal with essential endomorphisms of $(V^\cdot,A)$ from
the very beginning, i.e. with those sequences $E^i \in \mathcal{L} (V^i)$ which
satisfy 
$
   E^{i+1} A^{i} = A^{i} E^{i}
$
modulo compact operators from $V^i$ to $V^{i+1})$.
However, these latter don't act naturally on the cohomology of reduced complexes.

\section{Pseudodifferential operators}
\label{s.Pseudo}
\setcounter{equation}{0}

In this section we sketch a calculus of pseudodifferential operators on compact closed manifolds. 
We will follow the definitions of Shubin's classical book \cite{Shub87}. 

Let $U$ be an open set in $\R^n$ and $m$ a real number.
Denote by
   $\mathcal{S}^m (U \times \R^n)$ 
the space of all smooth functions $a$ on $U \times \R^n$ with the property that, 
for each 
   $\alpha, \beta \in \mathbb{Z}_{\geq 0}^n$ and 
   compact set $K \subset U$, 
there exists a constant $c_{\alpha,\beta,K} > 0$ satisfying
$$
   |\partial_x^{\alpha} \partial_{\xi}^{\beta} a(x,\xi)|
 \leq c_{\alpha,\beta,K}\, (1 + |\xi|^2)^{\frac{\scriptstyle m-|\beta|}
                                               {\scriptstyle 2}} 
$$ 
for all $(x,\xi) \in K \times \R^n$. 
The elements of $\mathcal{S}^m (U \times \R^n)$ are called symbols and those of
$$
   \mathcal{S}^{-\infty} (U \times \R^n)
 = \bigcap_m \mathcal{S}^m (U \times \R^n)
$$ 
smoothing symbols. 

To any symbol 
   $a \in \mathcal{S}^m (U \times \R^n)$ 
we assign the canonical pseudodifferential operator $A = a (x,D)$ by
$$
   Au\, (x) = \mathcal{F}^{-1}_{\xi \mapsto x}
              a (x,\xi)
              \mathcal{F}_{x \mapsto \xi} u  
$$
for $u \in C^{\infty}_{\mathrm{comp}} (U)$, where
   $\mathcal{F} u$ is the Fourier transform of $u$.  
Note that $A$ maps $C^{\infty}_{\mathrm{comp}} (U)$ continuously into 
                   $C^{\infty} (U)$.
The function $\sigma (A) := a$ is called the symbol of $A$.

We now want to consider classical pseudodifferential operators. 
They form an important subclass of canonical pseudodifferential operators which is closed under basic operations. 
Classical pseudodifferential operators were introduced in 1965 by 
   J.~J. Kohn and 
   L. Nirenberg
who reinforced the theory of S.~G. Michlin, A.~P. Calderon, etc.
The main property of this class is the existence of a principal symbol.
More precisely, a symbol 
   $a \in \mathcal{S}^m (U \times \R^n)$ 
is said to be classical (or multihomogeneous) if there is a sequence 
   $\{ a_{m-j} \}_{j = 0, 1, \ldots}$ 
of functions 
   $a_{m-j} \in C^{\infty} (U \times (\R^n \setminus \{ 0 \}))$
positively homogeneous of degree $m-j$ in $\xi$, such that 
$$
   a - \chi \sum_{j=0}^N a_{m-j} 
 \in \mathcal{S}^{m-N-1} (U \times \R^n)
$$
for all $N = 0, 1, \ldots$, where 
   $\chi \in C^{\infty} (\R^n)$ 
is a cut-off function with respect to $\xi = 0$.
Obviously, all the components $a_{m-j}$ are uniquely determined by $a$.
A canonical pseudodifferential operator $A$ on $U$ is called classical if its symbol $\sigma (A)$ is classical.
The set of all classical pseudodifferential operators of degree $m$ on $U$ is denoted by $\PsDO_{\mathrm{cl}}^m (U)$.      
The component $\sigma^m (A) := a_m$ is called the principal symbol of $A$.

\begin{example}
{\em
Any linear partial differential operator $A$ of order $m$ on $U$ has the form 
$$
 A (x,D)
  := \sum_{|\alpha| \leq m} A_{\alpha} (x) D^{\alpha},
$$  
where $A_{\alpha} \in C^{\infty} (U)$.
This is a classical pseudodifferential operator with symbol 
   $\sigma (A) (x,\xi) = A (x,\xi)$.
The principal symbol of $A$ is
$$
   \sigma^m (A) (x,\xi)
 = \sum_{|\alpha|= m} A_{\alpha} (x) \xi^{\alpha}.
$$  
}
\end{example}

Canonical pseudodifferential operators on open sets of $\R^n$ glue together to give rise to pseudodifferential operators on sections of vector bundles over a
smooth manifold $X$ of dimension $n$.
Locally sections of a vector bundle $E$ of rank $k$ are functions with values in $\C^k$.
Canonical pseudodifferential operators mapping functions with values in $\C^k$ 
                                            to functions with values in $\C^l$ are simply $(l \times k)\,$-matrices 
\begin{equation}
\label{eq.psdo}
   \Big( a_{i j} (x,D) \Big)_{i = 1, \ldots, l \atop
                              j = 1, \ldots, k}
\end{equation}
of canonical pseudodifferential operators on scalar-valued functions. 
The notions of multihomogeneity and principal symbol are extended to 
   (\ref{eq.psdo})
in entry wise unless a sophisticated approach is elaborated.
Given vector bundles $E$ and 
                     $F$ 
of ranks $k$ and 
         $l$ 
over $X$, by a pseudodifferential operator mapping sections of $E$ to those of $F$ is meant any map 
   $A : C^{\infty}_{\mathrm{comp}} (X,E) \to C^{\infty} (X,F)$
which has form (\ref{eq.psdo}) in any coordinate patch $U$ in $X$ over which both $E$ and $F$ are trivial,
   for any choice of local coordinates and trivializations. 
The space of all classical pseudodifferential operators of order $m$ on $X$ mapping sections of $E$ to those of
                    $F$
is denoted by $\PsDO^m_{\cl}(X; E,F)$. 
For $A \in \PsDO_{\cl}^m (X; E,F)$, the principal symbol $\sigma^m (A)$ proves to be a well-defined homomorphism of induced bundles
   $\pi^{\ast} E \to \pi^{\ast} F$
over $T^\ast X \setminus \{ 0 \}$, where
   $\pi\!: (T^{\ast} X \setminus \{ 0 \}) \to X$ is the natural projection.
Locally the principal symbol is a family of linear maps
   $\sigma^m (A) (x,\xi) \! : E_x \to F_x$
parametrized by $\xi \in \R^n \setminus \{ 0 \}$, where $E_x$ and
                                                        $F_x$
stand for the fibers of $E$ and $F$ over a point $x \in X$.
It is positively homogeneous of order $m$ in $\xi$. 

From now on we assume that $X$ is a compact closed smooth manifold of dimension
$n$ (e.g. $X$ is a sphere in $\R^{n+1}$).
In this case pseudodifferential operators on $X$ can be composed with each other thus giving rise to the simplest operator algebra 
$$
   \bigcup_{m \in \R} \PsDO^m_{\cl} (X).
$$      

\begin{theorem}
\label{t.composition}
Suppose that
   $A \in \PsDO^m_{\cl} (X;E,F)$ and 
   $B \in \PsDO^n_{\cl} (X;F,G)$. 
Then $BA \in \PsDO^{m+n}_{\cl} (X;E,G)$ and 
$
   \sigma^{m+n} (BA) = \sigma^n (B) \sigma^m (A)
$
holds.
\end{theorem}

Let $A \in \PsDO_{\cl}^m (X; E,F)$. 
On endowing the bundles $E$ and
                        $F$
by Riemannian metrics we introduce the formal adjoint 
   $A^{\adj} \! : C^{\infty} (X,F) \to C^{\infty} (X,E)$ 
by requiring
$$
   (Au,g)_{L^2 (X,F)} = (u, A^\adj g)_{L^2 (X,E)}
$$
for all $u \in C^{\infty} (X,E)$ and
        $g \in C^{\infty} (X,F)$.
A priori it is by no means clear if any $A^\adj$ exists.

\begin{theorem}
\label{t.adjoint}
Each operator
   $A \in \PsDO^m_{\cl} (X;E,F)$ 
possesses a formal adjoint $A^\adj \in \PsDO^m_{\cl} (X;F,E)$ and 
$
   \sigma^{m} (A^{\adj}) = (\sigma^m (A))^\ast
$
holds.
\end{theorem}

An operator $A \in\PsDO_{\cl}^m (X; E,F)$ is said to be elliptic if 
   $\sigma^m (A) (x,\xi)\!: E_x \to F_x$ 
is invertible for all $x \in X$ and
                      $\xi \in T^\ast_x X \setminus \{ 0 \}$. 
 
\begin{theorem}
\label{S.Parametrix Ex}
For each elliptic operator 
   $A \in \PsDO_{\cl}^m (X;E,F)$
there is an operator $P \in \PsDO^{-m}_{\cl} (X;F,E)$, such that
\begin{equation}
\label{eq.parametrix}
   \begin{array}{rcl}
   \Id_{E} - PA 
 & \in 
 & \PsDO^{-\infty} (X;E),
\\
   \Id_{F} - AP
 & \in 
 & \PsDO^{-\infty} (X;F).
   \end{array}
\end{equation}
\end{theorem}

Any operator $P \in \PsDO^{-m}_{\cl} (X;F,E)$ satisfying equalities (\ref{eq.parametrix}) is said to be a formal parametrix of $A$. 
Theorem \ref{S.Parametrix Ex} gains in interest if we realize that any formal parametrix is actually a parametrix in the sense of Banach spaces.
In order to show this we first mention the so-called property of spectral invariance of the algebra of pseudodifferential operators on a compact closed smooth manifold.

\begin{theorem}
\label{t.spectral} 
Let $A \in \PsDO^m_{\cl} (X;E,F)$ be elliptic and invertible on smooth sections. Then $A^{-1} \in \PsDO^{-m}_{\cl} (X;F,E)$.  
\end{theorem}

On combining duality arguments and Theorem \ref{t.adjoint} one readily sees that
any pseudodifferential operator
   $A \in \PsDO^m_{\cl} (X;E,F)$ 
extends by continuity to a continuous mapping  
   $A \!: (C^\infty (X,E))' \to (C^\infty (X,F))'$
of spaces of distribution sections.
Since $X$ is compact each distribution on $X$ is of finite order.
Hence it follows that the space $(C^\infty (X,E))'$ is exhausted by the scale of Sobolev spaces $H^s (X,E)$ with $s \in \R$, 
   and similarly for $(C^\infty (X,F))'$.
Pick a formally selfadjoint operator
   $\iL_E \in \PsDO^{2}_{\cl} (X;E)$ 
which is nonnegative and invertible on smooth sections of $E$
(e.g. $\iL_E = \partial_E^{\adj} \partial_E + \Id_E$ where $\partial_E$ is a 
 connection on $E$).
For $s \in \R$, the Sobolev space $H^s (X,E)$ is defined to consist of all  
$
   u \in (C^\infty (X,E))',
$
such that $\iL_E^{s/2}u \in L^2 (X,E)$.
This is a Hilbert space with scalar product
$$
   (u,v)_{H^s (X,E)}
 :=
   (\iL_E^{s/2} u, \iL_E^{s/2} v)_{L^2 (X,E)}
$$
for $u, v \in H^s (X,E)$. 

\begin{theorem}
\label{S.Fortsetzung H} 
For each $s \in \R$, any operator 
   $A \in \PsDO^m_{\cl} (X;E,F)$ 
maps $H^s (X,E)$ continuously into $H^{s-m} (X,F)$.
\end{theorem}

Combining this theorem with Rellich's theorem on compact embeddings of Sobolev
spaces we obtain the first defining property of the principal symbol mapping.
Together with the multiplicativity property of Theorem \ref{t.composition} this
allows one to identify the principal symbol mapping with the functor $\phi_{\Sigma}$ of Section \ref{Fred}.

\begin{theorem}
\label{t.Rellich} 
Let $A \in \PsDO^m_{\cl} (X;E,F)$ satisfy the condition $\sigma^m (A) = 0$. 
Then the extension $A\!: H^s (X,E) \to H^{s-m} (X,F)$ is a compact operator for all $s \in \R$. 
\end{theorem}

Suppose
   $A \in \PsDO^m_{\cl} (X; E,F)$ 
is an elliptic operator invertible on smooth sections. 
Then the extension 
   $A\!: H^s (X,E) \to H^{s-m} (X,F)$ 
is invertible, too. 
To show this we use the fact that the inverse $A^{-1}$ on smooth sections is actually a pseudodifferential operator in $\PsDO^{-m}_{\cl} (X;F,E)$, which is due to Theorem \ref{t.spectral}. 
By assumption, 
\begin{equation}
\label{eq.inverse}
\begin{array}{rcl}
   A^{-1}\, A u
 & = 
 & u,
\\
   A\, A^{-1} f
 & =
 & f
\end{array}
\end{equation}
hold for all $u \in C^\infty (X,E)$ and
             $f \in C^\infty (X,F)$.
Let $u \in H^s (X,E)$. 
Since $C^\infty (X,E)$ is dense in $H^s (X,E)$, there is a sequence $\{ u_j \}$ in $C^\infty (X,E)$ which converges to $u$ in $H^s (X,E)$. 
By Theorem \ref{S.Fortsetzung H}, $A$ maps $H^s (X,E)$ continuously to
                                           $H^{s-m} (X,E)$
and $A^{-1}$ maps $H^{s-m} (X,F)$ continuously to
                  $H^{s} (X,E)$.
Hence it follows that
$$
   A^{-1}\, A u = \lim A^{-1}\, A u_j = \lim u_j = u,
$$
i.e. the first equality in (\ref{eq.inverse}) is valid for all 
   $u \in H^s (X,E)$.
In the same one sees that the second equality in (\ref{eq.inverse}) holds true for all  
   $f \in H^{s-m} (X,F)$,
as desired.

\begin{theorem}
\label{S.Kern} 
Suppose that $A \in \PsDO^m_{\cl} (X;E,F)$ is elliptic.
Then the operator 
$
   A \! : H^s (X,E) \to H^{s-m} (X,F)
$
is Fredholm and its null-space belongs to $C^\infty (X,E)$.  
\end{theorem}

By the above, any $A \in \PsDO^m_{\cl} (X;E,F)$ gives rise to a continuous linear operator 
   $A \! : H^s (X,E) \to H^{s-m} (X,F)$,
for each fixed $s \in \R$.
Denote by 
$
   A^{\ast} \! : H^{s-m} (X,F) \to H^s (X,E)
$ 
the adjoint of this operator in the sense of Hilbert spaces, the adjoint  depending obviously on $s$.
We now specify this operator.

\begin{theorem}
\label{S.adjoint} 
If $A \in \PsDO^m_{\cl} (X;E,F)$ then    
$
   A^\ast  
 = \iL_E^{-s} A^\adj \iL_F^{s-m}.
$  
In particular, $A^\ast \in \PsDO^{-m}_{\cl} (X;F,E)$.
\end{theorem}

{\bf Proof.}
Since 
   $\iL_E^{-s} \in \PsDO^{-2s}_{\cl} (X;E)$ and
   $\iL_F^{s-m} \in \PsDO^{2 (s-m)}_{\cl} (X;F)$,
it suffices to establish the equality
$
   A^\ast  
 = \iL_E^{-s} A^\adj \iL_F^{s-m}
$
only.
By the definition of Hilbert adjoint we obtain 
\begin{eqnarray*}
   (Au,g)_{H^{s-m} (X,F)}   
 & = &
   (\iL_F^{(s-m)/2} Au, \iL_F^{(s-m)/2} g)_{L^2 (X,F)}
\\
 & = &
   (Au, \iL_F^{s-m} g)_{L^2(X,F)}
\\
 & = &
	(u, A^\adj \iL_F^{s-m} g)_{L^2 (X,E)}
\\
 & = &
	(\iL_E^{s/2} u, \iL_E^{s/2} \iL_E^{-s} A^\adj \iL_F^{s-m} g)_{L^2 (X,E)}
\\
 & = &
	(u, \iL_E^{-s} A^\adj \iL_F^{s-m} g)_{H^s (X,E)}
\end{eqnarray*}
for all $u \in C^\infty (X,E)$ and 
        $g \in C^\infty (X,F)$.
So the assertion follows by a familiar density argument.
\hfill $\square$

\section{Index theory for elliptic complexes}
\label{s.Index}
\setcounter{equation}{0}

Let $X$ be a compact closed smooth manifold of dimension $n$.
Consider a complex of pseudodifferential operators  
   $A^i \in \PsDO_{\cl}^{m_i} (X;F^i,F^{i+1})$
over $X$, i.e. a sequence 
$$ 
   C^\infty (X,F^{\cdot})\!:\,\,  
   0 
 \rightarrow 
   C^\infty (X,F^{0}) 
 \stackrel{A^0}{\rightarrow}  
   C^\infty (X,F^{1}) 
 \stackrel{A^{1}}{\rightarrow}   
   \ldots
 \stackrel{A^{N-1}}{\rightarrow}  
   C^\infty (X,F^{N}) 
 \rightarrow 
   0
$$
with the property
$
A^{i+1} A^i = 0.
$ 
To this complex we may assign the complex of principal symbols 
$$ 
   \pi^{\ast} F^{\cdot}\!:\,\,  
   \,0\, 
 \rightarrow
   \,\pi^{\ast} F^{0}\, 
 \stackrel{\sigma^{m_0} (A^0)}{\rightarrow}  
   \,\pi^{\ast} F^{1}\, 
 \stackrel{\sigma^{m_1} (A^1)}{\rightarrow}    
   \,\ldots\,
 \stackrel{\sigma^{m_{N-1}} (A^{N-1})}{\rightarrow}   
   \,\pi^{\ast} F^{N}\, 
 \rightarrow
   \,0.
$$

The complex $C^\infty (X,F^{\cdot})$ is called elliptic if the symbol complex $\pi^\ast F^{\cdot}$ is exact (away from the zero section of $T^\ast X$). 
In Section \ref{s.Elliptisch} we show that the cohomology of elliptic complexes is finite dimensional, a result going back at least as far as \cite{AtiyBott67}.

\begin{example}
{\em 
The de Rham complex $\iO^\cdot (X)$ of $X$ is elliptic.
For the complete bibliography we refer the reader to \cite{Tark90}.
}
\end{example}

In order to state the important Atiyah-Singer index formula we need some basic facts from differential geometry.  
Let $F$ be a smooth vector bundle over $X$. 
By a connection on $F$ is meant a first order differential operator      
   $\partial \!: C^{\infty} (X,F) \to \iO^1 (X,F)$ 
satisfying the Leibniz rule
$
   \partial (f u) = df\, u + f\, \partial u
$
for all $u \in C^\infty (X,F)$ and 
        $f \in C^\infty (X)$. 
The Leibniz rule allows one to extend the connection to differential forms of degree $q$ with coefficients in $F$ on $X$. 
We thus arrive at the sequence 
$$ 
   \iO^\cdot (X,F) :\,\,
   0\, 
 \rightarrow
   \,C^\infty (X,F)\, 
 \stackrel{\partial}{\rightarrow}  
   \,\iO^1 (X,F)\,
 \stackrel{\partial}{\rightarrow}   
   \,\ldots\,
 \stackrel{\partial}{\rightarrow}   
   \,\iO^n (X,F)\,
 \rightarrow
   \,0.
$$
The operator $\iO = \partial^2$ is a differential operator of order $0$ and is called the curvature. 
More precisely this is a matrix with items in differential forms of degree $2$
on $X$. 

If $f (z)$ is an analytic function in a neighborhood of $z = 0$ then it expands in powers of $z$.
On substituting $z = \iO$ we define $f (\iO)$, provided that the power series converges.
This is the case indeed, for $\iO^k$ vanishes, if $2k$ exceeds the dimension of $X$, and so the power series breaks.  
The characteristic classes of the bundle $F$ are defined by using a curvature of $F$.
They are actually independent modulo cohomology on the particular choice of curvature $\partial$.  
For example,
$$
\begin{array}{rcl}
   \mathrm{ch}\, (F)
 & =
 & \mathrm{tr}\, \exp \omega,
\\
   \mathcal{T}\, (F)
 & =
 & \displaystyle \det \Big( \frac{\omega}{1 - \exp (- \omega)} \Big)
\end{array}
$$
with 
$$
   \omega = - \frac{\iO}{2 \pi \imath}
$$
are the Chern character $\mathrm{ch} (F)$ and the Todd class $\mathcal{T} (F)$ of $F$, see for instance \cite{Pala65}, 
                         \cite{Well80}.  

Let 
   $A \in \PsDO^m_{\cl} (X;E,F)$ be an elliptic pseudodifferential operator on 
   $X$ 
and
   $a : \pi^\ast E \to \pi^\ast F$ the principal symbol of $A$.  
Then, $a$ induces an isomorphism of vector bundles  
   $\pi^\ast E$ and
   $\pi^\ast F$
over $T^\ast X \setminus \{ 0 \}$.
In this way $A$ determines an element $d (A)$ of the functor with compact support $K^{\mathrm{comp}} (T^\ast X)$.
Namely $d (A)$ is given by the virtual bundle with compact support 
   $\{ E, F, a \}$.
It is called the differential bundle of the elliptic operator.
The Chern character extends to the virtual bundle by
   $\mathrm{ch}\, (d (A)) := \mathrm{ch} (E) - \mathrm{ch} (F)$
(see \cite{Pala65} for a thorough treatment).

Assume $X$ is orientable. 
One can define the topological index $\mathrm{ind_{\mathrm{top}}}(A)$ of $A$ by 
$$
   \mathrm{ind}_{\mathrm{top}} (A)
 = \int_{T^{\ast} X} 
   \mathrm{ch}\, (d (A))\, \mathcal{T}\, (T_\C X),
$$  
where the orientation of the manifold $T^{\ast} X$ is given by the differential form
   $d \xi_1 \wedge dx^1 \wedge \ldots \wedge d \xi_n \wedge dx^n$.

\begin{theorem} [Atiyah-Singer index formula]
\label{t.index} 
Let 
   $X$ be oriented and
   $A$ be an elliptic operator on $X$. 
Then the analytical and the topological index coincide,
i.e. 
$  
   \mathrm{ind}\, (A)
 = \mathrm{ind}_{\mathrm{top}} (A).
$ 
\end{theorem}

{\bf Proof}
See \cite{Pala65}.
\hfill $\square$

\medskip

If $C^\infty (X, F^{\cdot})$ is an elliptic complex over $X$, we split 
   $V = \oplus C^\infty (X,F^{i})$ 
into the sum
$$ 
   V = V^{\mathrm{even}} \oplus V^{\mathrm{odd}}
$$
with 
   $V^{\mathrm{even}} = \oplus C^\infty (X,F^{2i})$ and 
   $V^{\mathrm{odd}} =\oplus C^\infty (X,F^{2i+1})$ 
and consider the block operator 
$$
   (A \oplus A^\adj)_e :\, 
   C^\infty (X,\oplus F^{2i}) \to C^\infty (X,\oplus F^{2i+1}) 
$$ 
given by 
$$
   \left( \begin{array}{ccccc}
          A^0 	  & A^1{}^\adj &	0 			   &	 	0 		  & \ldots 
\\
          0		  & A^2 		   & A^3{}^\adj 	& 		0 		  & \ldots
\\
          0 	  & 0		  	   & A^4 	      &  A^5{}^\adj & \ldots
\\
			 \ldots & \ldots		& \ldots	      &  \ldots     & \ldots
          \end{array}
   \right).
$$
Then $(A \oplus A^\adj)_e$ is elliptic and for the Euler characteristic of the complex we find
$$
   \chi (C^\infty (X,F^{\cdot}))
 = \int_{T^{\ast} X} 
   \mathrm{ch}\, (d (A \oplus A^\adj)_e)\, \mathcal{T}\, (T_\C X).
$$

It should be noted that for elliptic complexes of pseudodifferential operators of different order the ellipticity of $(A \oplus A^\adj)_e$ is understood in the sense of Douglis-Nirenberg.

\section{Elliptic quasicomplexes}
\label{s.Elliptisch}
\setcounter{equation}{0}

Having disposed of the preliminary material on elliptic complex, we turn to the
concept of elliptic quasicomplexes.
As above, $X$ stands for a compact closed $C^\infty$ manifold of dimension $n$.

By a quasicomplex of pseudodifferential operators on $X$ is meant any sequence of the form 
\begin{equation}
\label{eq.quasicomplex} 
   C^\infty (X,F^{\cdot})\!:\,\,  
   0 
 \rightarrow 
   C^\infty (X,F^{0}) 
 \stackrel{A^0}{\rightarrow}  
   C^\infty (X,F^{1}) 
 \stackrel{A^{1}}{\rightarrow}   
   \ldots
 \stackrel{A^{N-1}}{\rightarrow}  
   C^\infty (X,F^{N}) 
 \rightarrow 
   0
\end{equation}
with
   $A^i \in \PsDO_{\cl}^{m_i} (X;F^i,F^{i+1})$ 
satisfying
$
   A^{i+1} A^i \in \PsDO^{-\infty} (X; F^i,F^{i+2}).
$
In other words, the curvature of $C^\infty (X,F^\cdot)$ is a smoothing operator in the operator algebra under study.
To any quasicomplex one assigns the sequence of principal symbols 
$$ 
   \pi^{\ast} F^{\cdot}\!:\,\,  
   \,0\, 
 \rightarrow
   \,\pi^{\ast} F^{0}\, 
 \stackrel{\sigma^{m_0} (A^0)}{\rightarrow}  
   \,\pi^{\ast} F^{1}\, 
 \stackrel{\sigma^{m_1} (A^1)}{\rightarrow}    
   \,\ldots\,
 \stackrel{\sigma^{m_{N-1}} (A^{N-1})}{\rightarrow}   
   \,\pi^{\ast} F^{N}\, 
 \rightarrow
   \,0
$$
which is actually a complex of bundle homomorphisms, for
\begin{eqnarray*}
  \sigma^{m_{i+1}}(A^{i+1})\, \sigma^{m_i}(A^{i})  
 & = &
  \sigma^{m_{i+1}+m_i} (A^{i+1} A^{i}) 
\\
 & = &
   0.
\end{eqnarray*}

\begin{definition} 
\label{d.elliptic}
A quasicomplex $C^\infty (X,F^{\cdot})$ is called elliptic if its symbol complex $\pi^\ast F$ is exact (away from the zero section of $T^\ast X$). 
\end{definition}

If the orders of all operators $A^i$ in sequence (\ref{eq.quasicomplex}) are the same, 
   that is $m_i = m$ for all $i = 0, 1, \ldots, N-1$,
then one introduces the formal Laplacians of the quasicomplex 
   $C^\infty (X,F^{\cdot})$ 
by  
$
   L^i = A^{i-1} A^{i-1}{}^\adj + A^i{}^\adj A^i.
$
The quasicomplex $C^\infty (X,F^\cdot)$ is elliptic if and only if all the formal Laplacians $L^i$ are elliptic, as is easy to check.
Indeed, the principal symbols of $L^i$ just amount to the Laplacians of the
symbol complex $\pi^\ast F^{\ast}$.
However, this is no longer true if the orders $m_i$ are different, in which case formal Laplacians fail to be relevant to the study.

To get rid of this irrelevance we use a familiar construction with order reduction isomorphisms. 
More precisely, for a fixed $s \in \R$, we choose invertible operators 
$$
   R_i \in \PsDO_{\cl}^{s-(m_0 + \ldots + m_{i-1})} (X; F^i)
$$ 
($m_{-1} = 0$) and set 
$
   \tilde{A}{}^i = R_{i+1} A^i R_i^{-1}.
$
Then 
   $\tilde{A}{}^i \in \PsDO_{\cl}^{0} (X;F^i,F^{i+1})$ 
holds and we obtain a quasicomplex 
$$
   (C^\infty (X,F^{\cdot}), \tilde{A}) \! : \,\,  
   0 
 \rightarrow 
   C^\infty (X,F^{0}) 
 \stackrel{\tilde{A}{}^0}{\rightarrow}  
   C^\infty (X,F^{1}) 
 \stackrel{\tilde{A}{}^1}{\rightarrow}   
   \ldots
 \stackrel{\tilde{A}{}^{N-1}}{\rightarrow}  
   C^\infty (X,F^{N}) 
 \rightarrow 
   0 
$$
of operators of degree $0$.

\begin{theorem}
\label{S.ellipt gdw}
A quasicomplex $C^\infty (X,F^{\cdot})$ is elliptic if and only if the reduced quasicomplex $(C^\infty (X,F^{\cdot}),\tilde{A})$ is elliptic.
\end{theorem}

{\bf Proof.}
Set 
\begin{equation}
\label{eq.reduction}
   s_i = s - (m_0 + \ldots + m_{i-1})
\end{equation} 
for $i = 0, 1, \ldots, N$, so that $s_0 := s$ and 
                                   $s_{i+1} = s_i - m_i$.
Since 
   $R_i \in \PsDO_{\cl}^{s_i} (X; F^i)$ 
is invertible, the inverse operator is available in 
   $\PsDO_{\cl}^{-s_i} (X; F^i)$, 
which is a consequence of spectral invariance. 
Moreover, we get 
$
   \sigma^{-s_i} (R_{i}^{-1})
 = (\sigma^{s_i} (R_{i}))^{-1}
$
whence
$$
   \sigma^{0} (\tilde{A}{}^{i}) 
 = \sigma^{s_{i+1}} (R_{i+1})\,
   \sigma^{m_i} (A^{i})\, 
   (\sigma^{s_i} (R_{i}))^{-1}.
$$
With this equality the assertion follows by a straightforward computation, as desired.
\hfill $\square$

\medskip

We may extend the quasicomplex $C^\infty (X,F^{\cdot})$ to a quasicomplex with Sobo\-lev spaces, i.e.
$$ 
   H^{s_{\cdot}} (X,F^{\cdot})\!:\,\,  
   \,0\, 
 \rightarrow 
   \,H^{s_0} (X,F^{0})\, 
 \stackrel{A^0}{\rightarrow}  
   \,H^{s_1} (X,F^{1})\, 
 \stackrel{A^{1}}{\rightarrow}   
   \,\ldots\,
 \stackrel{A^{N-1}}{\rightarrow}  
   \,H^{s_N} (X,F^{N}) 
 \rightarrow 
   \,0
$$ 
where $s_i$ are given by (\ref{eq.reduction}).      
From Theorem \ref{t.Rellich} it follows that 
   $H^{s_{\cdot}} (X,F^{\cdot})$ 
is a quasicomplex in the context of Hilbert spaces.
Theorem \ref{S.adjoint} shows readily that the Laplacians  
$
   \iD^i := A^{i-1} A^{i-1}{}^\ast + A^i{}^\ast A^i
$
of this quasicomplex belong to $\PsDO_{\cl}^{0} (X;F^i)$.

\begin{theorem}
\label{S.Laplace gdw}
A quasicomplex $C^\infty (X,F^{\cdot})$ is elliptic if and only if the Lap\-lacians  $\iD^i$ are elliptic.
\end{theorem}

{\bf Proof.}
Let
$
   \iL_{F^i} \in \PsDO_{\cl}^{2} (X;F^i)
$
be those invertible operators which define the norm in $H^s (X,F^i)$.
Set 
$$
   R_{i} := \Lambda^{s_i/2}
$$ 
and
$
   \tilde{A}{}^i := R_{i+1} A^i R_i^{-1}.
$
By definition, the diagram
$$
\begin{array}{ccccccccccc}
   0
 & \rightarrow
 & H^{s_0} (X,F^0)
 & \stackrel{A^{0}}\rightarrow
 & H^{s_1} (X,F^1)
 & \stackrel{A^{1}}\rightarrow
 & \ldots
 & \stackrel{A^{N-1}}\rightarrow
 & H^{s_N} (X,F^N)
 & \rightarrow
 & 0
\\
 & 
 & \mbox{ }\downarrow R_0 
 & 
 & \mbox{ }\downarrow R_1
 & 
 & 
 & 
 & \mbox{ }\downarrow R_N
 & 
 & 
\\
   0
 & \rightarrow
 & L^2 (X,F^0)
 & \stackrel{\tilde{A}{}^{0}}\rightarrow
 & L^2 (X,F^1)
 & \stackrel{\tilde{A}{}^{1}}\rightarrow
 & \ldots
 & \stackrel{\tilde{A}{}^{N-1}}\rightarrow
 & L^2 (X,F^N)
 & \rightarrow
 & 0
\end{array}
$$
is commutative.
With
\begin{eqnarray*}
   A^i{}^\ast 
 & = &
   \iL_{F^i}^{-s_i}\, A^{i}{}^\adj\, \iL_{F^{i+1}}^{s_{i+1}}
\\
 & = &
   R_i^{-2}\, A^{i}{}^\adj\, R_{i+1}^{2} 
\end{eqnarray*}
and
$
   A^{i}{}^\adj 
 = R_i\, \tilde{A}{}^{i}{}^\adj\, R_{i+1}^{-1} 
$
we find
\begin{eqnarray*}
   \iD^i 
 & = &
   ( R_{i}^{-1} \tilde{A}{}^{i-1} R_{i-1} )\,  
   ( R_{i-1}^{-1}\, \tilde{A}{}^{i-1}{}^\adj\, R_{i} ) 
 + ( R_i^{-1}\, \tilde{A}{}^{i}{}^\adj\, R_{i+1} )\, 
   ( R_{i+1}^{-1} \tilde{A}{}^{i} R_{i} )   
\\
 & = &
   R_{i}^{-1}\, \tilde{\iD}{}^i\, R_{i}, 
\end{eqnarray*}
where 
$
   \tilde{\iD}{}^i
 = \tilde{A}{}^{i-1} \tilde{A}{}^{i-1}{}^\adj 
 + \tilde{A}{}^{i}{}^\adj \tilde{A}{}^{i}
$ 
are the Laplacians of the reduced quasicomplex.
Since 
$$
 \sigma^{0} (\iD^{i}) 
 = (\sigma^{s_i} (R_{i}))^{-1}\,
   \sigma^{0} (\tilde{\iD}{}^{i})\, 
   \sigma^{s_i} (R_{i}),
$$
the pseudodifferential operator $\iD^i$ is elliptic if and only is so is
                                $\tilde{\iD}{}^{i}$.
Applying Theorem \ref{S.ellipt gdw} we get the desired assertion.
\hfill $\square$

By a formal parametrix of quasicomplex $C^\infty (X,F^{\cdot})$ is meant any sequence of pseudodifferential operators
   $P^{i} \in \PsDO^{-m_{i-1}} (X; F^{i},F^{i-1})$ 
satisfying the homotopy equations 
$$ 
   P^{i+1} A^{i} + A^{i-1} P^{i} = \Id_{F^{i}} - S^{i}
$$ 
with smoothing operators 
   $S^{i} \in \PsDO^{-\infty} (X; F^{i})$ 
for all $i = 0, 1, \ldots, N$.

\begin{theorem}
Each elliptic quasicomplex $C^\infty (X,F^{\cdot})$ possesses a formal parametrix.
\end{theorem}

{\bf Proof.} 
We consider the extended quasicomplex $H^{s_{\cdot}} (X,F^{\cdot})$.
By Theorem \ref{S.Laplace gdw}, the Laplacians $\iD^i$ are elliptic, and so we use Theorem \ref{S.Parametrix Ex} to find a formal parametrix 
   $G^i \in \PsDO_{\cl}^{0} (X;F^i)$ 
for $\iD^i$. 
On arguing as in the proof of Theorem \ref{S.Fred Quasi} we show that the family 
   $P^i = G^{i-1} A^{i-1}{}^\ast$ 
is a formal parametrix for $C^\infty (X,F^{\cdot})$, as desired. 
\hfill $\square$

\begin{theorem}
Assume that $C^\infty (X,F^{\cdot})$ is an elliptic quasicomplex. 
Then the extended quasicomplex $H^{s_\cdot} (X,F^{\cdot})$ is Fredholm.
\end{theorem}

{\bf Proof.}
Since $C^\infty (X,F^{\cdot})$ is elliptic, this quasicomplex possesses a formal parametrix. 
The extension of this parametrix to Sobolev spaces is a genuine parametrix of 
   $H^{s_\cdot} (X, F^{\cdot})$ 
in the sense of Hilbert spaces, for the smoothing operators are compact.
\hfill $\square$

\begin{theorem}
\label{t.reduction}
For each elliptic quasicomplex 
   $(C^\infty (X,F^{\cdot}),A)$  
there is an elliptic complex  
   $(C^\infty (X,F^{\cdot}),D)$, 
such that
   $D^i \in \PsDO_{\cl}^{m_i} (X; F^i,F^{i+1})$ and  
   $D^i = A^i$ modulo $\PsDO^{-\infty} (X; F^i,F^{i+1})$.
\end{theorem}

{\bf Proof.}
Since the differential $A$ of the quasicomplex is given by pseudodifferential
operators on $X$, we can turn to the quasicomplex 
   $(H^{s_{\cdot}} (X, F^{\cdot}),A)$
of Hilbert spaces, where $s \in \R$ is any fixed number.
Arguing as in the proof of Theorem 8.1 in \cite{KrupTarkTuom07} we modify the
quasicomplex $(H^{s_{\cdot}} (X, F^{\cdot}),A)$ into a complex 
   $(H^{s_{\cdot}} (X, F^{\cdot}),D)$
whose differential $D$ differs from $A$ by a smoothing operator.

We start from the end of the quasicomplex
$$
   \ldots
 \stackrel{A^{N-2}}{\rightarrow}
   H^{s_{N-1}} (X,F^{N-1})
 \stackrel{D^{N-1}}{\rightarrow}
   H^{s_N} (X,F^N)
 \rightarrow
   0,
$$
setting $D^{N-1} = A^{N-1}$.
Since $\sigma^{m_{N-1}} (A^{N-1})$ is surjective, it follows that the Laplacian
$
   \iD^{N}
 = D^{N-1} D^{N-1}{}^*
$
is a selfadjoint elliptic pseudodifferential operator of order $0$ in 
   $H^{s_N} (X,F^N)$.
By the Hodge theory for single operators, there is an operator
   $G^N \in \PsDO^{0} (X; F^N)$
satisfying
$$
   \Id_{F^N} 
 =
   H^N + \iD^N G^N
 =
   H^N + D^{N-1} P^N,
$$
where
   $H^N$ stands for the orthogonal projection onto the
   finite-dimensional space $\ker \iD^N = \ker D^{N-1}{}^\ast$
and
   $P^N = D^{N-1}{}^\ast G^N$.
We set
$$
   \iP^{N-1} = \Id_{F^{N-1}} - P^N D^{N-1}
$$
thus obtaining a pseudodifferential operator in $\PsDO^0 (X;F^{N-1})$.
We claim that
   $\iP^{N-1}$ is a projection onto $\ker D^{N-1}$.
Indeed,
   $\iP^{N-1} = \Id_{F^{N-1}}$ is valid on $\ker D^{N-1}$
and
\begin{eqnarray*}
   \iP^{N-1} \iP^{N-1}
 & = &
   (\Id_{F^{N-1}} - P^N D^{N-1}) (\Id_{F^{N-1}} - P^N D^{N-1})
\\
 & = &
   \Id_{F^{N-1}} 
 - 2\, P^N D^{N-1}
 + P^N (D^{N-1} P^N) D^{N-1}
\\
 & = &
   \Id_{F^{N-1}}
 - 2\, P^N D^{N-1}
 + P^N (\Id_{F^{N}} - H^N) D^{N-1}
\\
 & = &
   \iP^{N-1},
\end{eqnarray*}
for
$
   H^N D^{N-1}
 = (D^{N-1}{}^\ast H^N)^\ast
 = 0.
$

Next we set
   $D^{N-2} = \iP^{N-1} A^{N-2}$.
Then
   $D^{N-2} \in \PsDO^{m_{N-2}} (X;F^{N-2},F^{N-1})$
and
   $D^{N-1} D^{N-2} = 0$,
for $\iP^{N-1}$ is a projection onto $\ker D^{N-1}$.
Furthermore, we get
\begin{eqnarray*}
   D^{N-2}
 & = &
   A^{N-2}
 - P^N A^{N-1} A^{N-2}
\\
 & = &
   A^{N-2}
\end{eqnarray*}
modulo $\PsDO^{-\infty} (X;F^{N-2},F^{N-1})$, 
   for the composition $A^{N-1} A^{N-2}$ is a smoothing operator.

Consider now a slightly modified quasicomplex
$$
   \ldots
 \stackrel{A^{N-3}}{\rightarrow}
   H^{s_{N-2}} (X,F^{N-2})
 \stackrel{D^{N-2}}{\rightarrow}
   H^{s_{N-1}} (X,F^{N-1})
 \stackrel{D^{N-1}}{\rightarrow}
   H^{s_N} (X,F^N)
 \rightarrow
   0.
$$
Since 
   the symbol complex of the initial quasicomplex is exact 
and
   the operators $D^i$ and $A^i$ have the same principal symbol for 
   $i = N-2, N-1$, 
the Laplacian
$
   \iD^{N-1}
 = D^{N-2} D^{N-2}{}^\ast + D^{N-1}{}^\ast D^{N-1}$
is a selfadjoint elliptic operator of order $0$ on $H^{s_{N-1}} (X,F^{N-1})$.
Using the Hodge theory for complexes, we deduce that there is an operator
   $G^{N-1} \in \PsDO^{0} (X; F^{N-1})$,
such that
\begin{eqnarray*}
   \Id_{F^{N-1}}
 & = &
   H^{N-1}
 + D^{N-2} D^{N-2}{}^\ast G^{N-1}
 + D^{N-1}{}^\ast G^{N} D^{N-1}
\\
 & = &
   H^{N-1}
 + D^{N-2} P^{N-1}
 + P^{N} D^{N-1}
\end{eqnarray*}
where
   $H^{N-1}$ is the orthogonal projection onto the null-space of
   $\Delta^{N-1}$
which is
   $\ker D^{N-2}{}^\ast \cap \ker D^{N-1}$,
and
   $P^{N-1} = D^{N-2}{}^\ast G^{N-1}$.
Then, we claim that
$$
   \iP^{N-2} = \Id_{F^{N-1}} - P^{N-1} D^{N-2}
$$
is the orthogonal projection onto $\ker D^{N-2}$.
Indeed,
   $\iP^{N-2}$ is the identity operator on $\ker D^{N-2}$.
Moreover,
\begin{eqnarray*}
   (\iP^{N-2})^2
 & = &
   \iP^{N-2}
 - P^{N-1} D^{N-2}
 + P^{N-1} (D^{N-2} P^{N-1}) D^{N-2}
\\
 & = &
   \iP^{N-2}
 - P^{N-1} D^{N-2}
 + P^{N-1} (\Id_{F^{N-1}} - H^{N-1} - P^{N} D^{N-1}) D^{N-2}
\\
 & = &
   \iP^{N-2}
 - P^{N-1} H^{N-1} D^{N-2}
\\
 & = &
   \iP^{N-2},
\end{eqnarray*}
since
$
   H^{N-1} D^{N-2}
 = (D^{N-2}{}^\ast H^{N-1})^\ast
$
vanishes.
Introducing
$
   D^{N-3}
 = \iP^{N-2} A^{N-3}
$
we thus obtain
   $D^{N-2} D^{N-3} = 0$
and
\begin{eqnarray*}
   D^{N-3}
 & = &
   A^{N-3}
 - P^{N-1} D^{N-2} A^{N-3}
\\
 & = &
   A^{N-3}
\end{eqnarray*}
modulo $\PsDO^{-\infty} (X;F^{N-3},F^{N-2})$, for
$
   D^{N\!-\!2} A^{N\!-\!3} 
 = A^{N\!-\!2} A^{N\!-\!3} + (D^{N\!-\!2} - A^{N\!-\!2}) A^{N\!-\!3}
$ 
is a smoothing operator.

Continuing in this fashion, in a finite number of steps we obtain a complex of
operators
   $D^i \in \PsDO^{m_i} (X; F^i,F^{i+1})$,
such that
   $D^i - A^i \in \PsDO^{-\infty} (X; F^i,F^{i+1})$
for all $i = 0, 1, \ldots, N-1$.
\hfill $\square$

It is worth pointing out that the desired complex $(C^\infty (X,F^\cdot), D)$ 
is constructed within the pseudodifferential calculus on $X$.
I.e., $D^i$ are pseudodifferential operators even in the case if the initial sequences of symbols stem from differential operators.

Let $C^\infty (X,F^{\cdot})$ be an elliptic quasicomplex. 
We define the Euler characteristic of this quasicomplex by
$$
   \chi (C^\infty (X,F^{\cdot})) := \chi (H^{s_{\cdot}} (X,F^{\cdot}),D),
$$
where $(H^{s_{\cdot}} (X,F^{\cdot}),D)$ is any complex with the properties listed in Theorem \ref{t.reduction}.

\begin{lemma}
As defined above, the Euler characteristic is independent neither of $s \in \R$ nor of the choice of the differential $D$ with $D^2 = 0$. 
\end{lemma}

{\bf Proof.}
It suffices to show that the cohomology of 
   $(H^{s_{\cdot}} (X,F^{\cdot}),D)$ 
is independent of the choice of $s$ and $D$. 
Since the Laplacians 
$
   \iD^{i} 
 =  D^{i-1} D^{i-1}{}^\ast + D^{i}{}^\ast D^{i}
$ 
are elliptic pseudodifferential operators, Theorem \ref{S.Kern} implies that the null-space of $\iD^i$ belongs to $C^{\infty} (X,F^{i})$, i.e. is independent of $s$.
Finally, the abstract Hodge theory yields
$$
   H^i (H^s (X,F^{\cdot}),D)
 \stackrel{\mathrm{top}}{\cong} 
   \ker \iD^i.
$$
The proof is completed by observing that the Euler characteristic of an elliptic complex 
   $(C^\infty (X,F^{\cdot}),D)$
is uniquely determined by the principal symbols of $D^i$, 
   see Theorem \ref{t.index}.
\hfill $\square$

Using the Atiyah-Singer index formula we are able to evaluate the index of an elliptic quasicomplex. 

\begin{theorem}   
Let $X$ be oriented and
   $C^\infty (X,F^{\cdot})$ 
be an elliptic quasicomplex. 
Then
$$
   \chi (C^\infty (X,F^{\cdot}))
 = \int_{T^{\ast} X} 
   \mathrm{ch}\, (d (A \oplus A^\adj)_e)\, \mathcal{T}\, (T_\C X)
$$
holds.
\end{theorem}

{\bf Proof.}
Since the analytical index does only depend on the principal symbols,
the assertion follows from the Atiyah-Singer index theorem for elliptic complexes.
\hfill $\square$

\medskip

One problem still unsolved concerns the Lefschetz fixed point formula for elliptic quasicomplexes.
Although particular cases have been well understood, no canonical concept of Lefschetz number is available in the context of quasicomplexes. 
One can also study quasicomplexes in the context of other operator algebras, e.g. on manifolds with singularities. 
In \cite{KrupTarkTuom07} this is done for compact manifolds with boundary.  
(Note that on manifolds with boundary one works with the Boutet de Monvel  
 algebra of pseudodifferential boundary value problems.) 
The next step is to consider elliptic quasicomplexes on compact manifolds with edges.

\vspace*{8mm}

\begin{scriptsize}

(Daniel Wallenta)
{\sc Universit\"{a}t Potsdam,
     Institut f\"{u}r Mathematik,
     Am Neuen Palais 10,
     14469 Potsdam,
     Germany}

E-mail address:
   wallenta@hotmail.de

\end{scriptsize}

\end{document}